\documentclass[12pt]{article}
\usepackage{graphicx}
\usepackage{amsmath}
\usepackage{epsf}
\usepackage{graphpap}

\numberwithin{equation}{section}

\def\be{\begin{equation}}
\def\ee{\end{equation}}

\date{February 02, 2008}
\begin{document}

\title{Use of Complex Lie Symmetries for Linearization of Systems of \\
Differential Equations - II: \\
Partial Differential Equations}
\author{S. Ali$^{a}$, F. M. Mahomed$^{b}$, Asghar Qadir$^{a}$ \\
$^{a}${\small Center For Advanced Mathematics and Physics, }\\
{\small National University of Sciences and Technology, }\\
{\small Rawalpindi, Pakistan.}\\
$^{b}${\small Centre for Differential Equations, Continuum }\\
{\small Mechanics and Applications, School of Computational }\\
{\small and Applied Mathematics, University of the}\\
{\small Witwatersrand, South Africa.}\\
}
\maketitle

\begin{abstract}
The linearization of complex ordinary differential equations is studied by
extending Lie's criteria for linearizability to complex functions of complex
variables. It is shown that the linearization of complex ordinary
differential equations implies the linearizability of systems of partial
differential equations corresponding to those complex ordinary differential
equations. The invertible complex transformations can be used to obtain
invertible real transformations that map a system of nonlinear partial
differential equations into a system of linear partial differential
equation. Explicit invariant criteria are given that provide procedures for
writing down the solutions of the linearized equations. A few nontrivial
examples are mentioned.
\end{abstract}

\section{\textbf{Introduction}}

\qquad The linearization i.e. use of point transformations to convert the
nonlinear equations to linear form of real ordinary differential equations
(RODEs) was first introduced by Lie \cite{lie}. He provided the
linearizability criteria for any scalar second order RODE to be point
transformable to a linear RODE via invertible maps of both independent and
dependent variables. He found that the RODE must be at most cubic in the
first derivative and that a particular overdetermined system of conditions
for the coefficents that must be satisfied (see e.g., [$10-13],[15]$).
Tresse \cite{tre} studied the linearization of a second order RODE by
looking at the two relative invariants of the equivalence group of
transformations, the vanishing of both of which gives the necassary and
sufficient conditions for linearization. These are equivalent to the Lie
conditions \cite{lea}.

Lie's approach has been extended in many ways. The linearization of higher
order RODEs and their systems has been discussed by various authors (see
e.g., \cite{lea, mah, sar, waf}). Further, linearization of PDEs is
discussed in \cite{blu, blu1, kum}. It is seen that for systems of RODEs,
Lie linearizability criteria involves is solving more than ten equations
simultaneously that become complicated and messy. Also, for PDEs and their
systems it is highly nontrivial to obtain the transformations that remove
nonlinearity in the equations.

We have extended the Lie approach of linearization of RODEs via complex
point transformations to the complex domain and obtained corresponding Lie
criteria for complex ordinary differential equations (CODEs). A CODE is a
differential equation of a single complex function of one complex variable.
The CODEs can be linearized in the same way as RODEs via complex point
transformations. A CODE yields a system of partial differential equation
(PDEs) after decomposing all the complex functions, variables and
derivatives into their real and imaginary parts. That also yield
Cauchy-Riemann equations (CR-equation) which are already linear thus yield no
difference for the linearization of systems of PDEs. The linearization of
systems of PDEs corresponding to a CODE follows directly from the
linearization of that CODE via real point transformations that are obtained
after decomposing the complex point transformation. This is an intrinsic
property of complex theory, which Penrose calls \textquotedblleft complex
magic\textquotedblright\ \cite{Pen}, that yields a non-trivial way of
linearizing PDEs and their systems.

The outline of the paper is as follows. In the second section we discuss
those CODEs that are equivalent, i.e. which can be transformed into one
another via invertible complex mappings. Lie compatibility conditions and
some other results for linearization of systems of PDEs corresponding to
CODEs are discussed in the third section. Examples are also discussed in the
same section. Finally, in the fourth section we have given a summary and
discussion regarding our complex Lie symmetry (CLS) analysis for
linearization.

\section{\textbf{Equivalent CODEs}}

\qquad Two CODEs are (locally) equivalent via an invertible complex
transformation if one can be transformed into the other by an invertible
complex transformation. A CODE is a system of four PDEs, which includes the
two linear Cauchy-Rimenann equations (CREs) in two unknown functions of two
indpendent variables. If a CODE is equivalent to another CODE via inveritble
complex transformations then the system of PDEs corresponding to that CODE
is also equivalent to the other system of PDEs. The invertible complex
transformation yields two real transformations, which can be used to
transform one system of PDEs to another system.

Every first-order CODE 
\begin{equation}
u^{\prime }(z)=w(z,u),  \tag{1}
\end{equation}%
can be transformed into the simplest one $u^{\prime }=0,$ via a complex
point tranformation \cite{mah1}. Similarly, every linear second-order CODE
can always be transformed into its simplest form via invertible complex
transformations, specifically into $u^{\prime \prime }(z)=0.$ For example,
the complexified Ricatti equation 
\begin{equation}
u^{\prime }(z)+u^{2}=0,  \tag{2}
\end{equation}%
is transformable to $U^{\prime }=0$ by means of 
\begin{equation}
Z=z,~U=(1/u)-z.  \tag{3}
\end{equation}%
Also, the complexified simple harmonic oscillator equation 
\begin{equation}
u^{\prime \prime }(z)+u=0,  \tag{4}
\end{equation}%
can be transformed into $U^{\prime \prime }=0$ via an invertible complex
transformation 
\begin{equation}
Z=\tan z,~U=u~\sec z.  \tag{5}
\end{equation}%
\newline
To decompose the above CODEs into corresponding systems of PDEs, write%
\begin{equation}
z =x+iy,~u(z)=f(x,y)+ig(x,y),  \tag{6} \\
\end{equation}%
\begin{equation}
w =w_{1}+iw_{2}.  \tag{7}
\end{equation}%
The system of PDEs corresponding to a first-order CODE is 
\begin{equation}
f_{x}+g_{y}=w_{1}(x,y,f,g),\text{ \ }g_{x}-f_{y}=w_{2}(x,y,f,g),  \tag{8}
\end{equation}%
which can be transformed into%
\begin{equation}
F_{x}+G_{y}=0,\text{ \ }G_{x}-F_{y}=0,  \tag{9}
\end{equation}%
by an invertible real transformation derived from the complex point
tranformation. Similarly, the system of PDEs corresponding to a linear
second-order CODE is%
\begin{equation}
f_{xx}-f_{yy}+2g_{xy}=w_{1}(x,y,f,g,h,l),  \tag{10}
\end{equation}%
\begin{equation}
g_{xx}-g_{yy}-2f_{xy}=w_{2}(x,y,f,g,h,l),  \tag{11}
\end{equation}%
where both $w_{1}$ and $w_{2}$ are two such real functions that do not give
rise to a nonlinear system. The above system can be transformed into%
\begin{equation}
F_{XX}-F_{YY}+2G_{XY}=0,  \tag{12}
\end{equation}%
\begin{equation}
G_{XX}-G_{YY}-2F_{XY}=0.  \tag{13}
\end{equation}%
via an invertible real transformation that is obtained from a complex
transformation. \newline
\newline
\textbf{Examples}\newline
\textbf{1.} The complexified Ricatti equation is equivalent to%
\begin{equation}
f_{x}+g_{y}=-f^{2}+g^{2},  \tag{14}
\end{equation}%
\begin{equation}
g_{x}-f_{y}=-2fg.  \tag{15}
\end{equation}%
If we apply the following transformation%
\begin{equation}
X=x,~Y=y,  \tag{16}
\end{equation}%
\begin{equation}
F=\frac{f}{f^{2}+g^{2}}-x,~G=\frac{-g}{f^{2}+g^{2}}-y,  \tag{17}
\end{equation}%
we get%
\begin{equation}
F_{X}+G_{Y}=0,  \tag{18}
\end{equation}%
\begin{equation}
G_{X}-F_{Y}=0.  \tag{19}
\end{equation}%
Notice that a simple complex transformation for the complexified Ricatti
equation yields a non-trivial real transformation (16) and (17) that
transforms the system of PDEs (14) and (15) into its simple analogue (18) and (19).
This remarkable feature of complex theory is of great
significance for the linearization of PDEs and their systems.\newline
\newline
\textbf{2.} The system%
\begin{equation}
f_{xx}-f_{yy}+2g_{xy}=-f,  \tag{20}
\end{equation}%
\begin{equation}
g_{xx}-g_{yy}-2f_{xy}=-g,  \tag{21}
\end{equation}%
can be transformed into%
\begin{equation}
F_{XX}-F_{YY}+2G_{XY}=0,  \tag{22}
\end{equation}%
\begin{equation}
G_{XX}-G_{YY}-2F_{XY}=0,  \tag{23}
\end{equation}%
via the real tranformations%
\begin{gather}
X=\frac{(1/2)\sin 2x}{\cos ^{2}x+\sinh ^{2}y},\quad Y=\frac{(1/2)\sinh 2y}{%
\cos ^{2}x+\sinh ^{2}y},  \tag{24} \\
F=\frac{f\cos x\cosh y-g\sin x\sinh y}{\cos ^{2}x+\sinh ^{2}y},\quad G=\frac{%
f\sin x\sin y+g\cos x\cosh y}{\cos ^{2}x+\sinh ^{2}y}.  \tag{25}
\end{gather}%
Thus, the system of PDEs corresponding to a CODE can easily be transformed
into its simplified (linear) form by invoking complex transformations. It
would have been very difficult to guess or calculate the real
transformations (24) and (25) that transform system $(20)$ and (21) into
system (22) and (23) without this formalism. It is an extraordinary
characteristic of complex magic.

Lie condititons in the complex domain can be derived by replacing the real
variable $x$ in r-CODE in \cite{Ali2} with a complex variable $z$. A scalar
second order CODE which is at most cubic in its first derivative%
\begin{equation}
u^{\prime \prime }(z)=A(z,u)u^{\prime 3}+B(z,u)u^{\prime 2}+C(z,u)u^{\prime
}+D(z,u),  \tag{26}
\end{equation}%
where $A$, $B,C$ and $D$ are complex valued functions is linearizable
according to Lie theorem. We do not re-state Theorem 1 of \cite{Ali2} as all
the conditions (1-9) remain the same with only the change that the real
variable $x$ is replaced by a complex variable $z,$ where, the constants in
(34) in \cite{Ali2} become complex constants for CODEs. In the subsequent
section we will decompose (26) into system of PDEs by using (6).

\section{\textbf{Lie Conditions for Systems of PDEs }}

\qquad According to Lie's theorem the necessary condition for a second order
CODE to be linearizable is that it is at most cubic in its first order
derivative implies the necessary conditions for system of PDEs corresponding
to that CODE i.e. the system must be at most cubic in first derivatives
together with certain constraints on coefficients. The real transformations
for linearization of a system of nonlinear PDEs can be obtained by
decomposing complex transformations that linearize a CODE. Now we derive Lie
conditions for systems of PDEs corresponding to a second order CODE by
decomposing results in Theorm 1 into real variables. The general form of a
system of PDEs corresponding to (26) is given by 
\begin{equation}
f_{xx}-f_{yy}+2g_{xy}=A_{1}(h^{3}-3hl^{2})-A_{2}(3h^{2}l-l^{3})+B_{1}(h^{2}-l^{2})-2B_{2}hl+C_{1}h-C_{2}L+D_{1},
\tag{27}
\end{equation}%
\begin{equation}
g_{xx}-g_{yy}-2f_{xy}=A_{1}(3h^{2}l-l^{3})+A_{2}(h^{3}-3hl^{2})+2B_{1}hl+B_{2}(h^{2}-l^{2})+C_{2}h+C_{1}L+D_{2}.
\tag{28}
\end{equation}%
where all the coefficients $A_{i},B_{i},C_{i}$ and $D_{i}$ are functions of $%
x,y,f,g.$\newline
\newline
\textbf{Theorem 1. }The following statements are equivalent.\newline
\textbf{1.} The above system of PDEs (27) and (28) is linearizable via real
transformations;\newline
\textbf{2.} The coefficients in (27) and (28) satisfy%
\begin{gather}
3A_{xx}^{1}-3A_{yy}^{1}+6A_{xy}^{2}+3C^{1}A_{x}^{1}+3C^{1}A_{y}^{2}-3A_{x}^{2}C^{2}+3C^{2}A_{y}^{1}-3A_{f}^{1}D^{1}-3D^{1}A_{g}^{2}+
\notag \\
3D^{2}A_{f}^{2}-3D^{2}A_{g}^{1}+3A^{1}C_{x}^{1}+3A^{1}C_{y}^{2}-3A^{2}C_{y}^{2}+3A^{2}C_{x}^{1}+C_{ff}^{1}-C_{gg}^{1}+2C_{fg}^{2}-
\notag \\
6A^{1}D_{f}^{1}-6A^{1}D_{g}^{2}+6A^{2}D_{f}^{2}-6A^{2}D_{g}^{1}+B^{1}C_{f}^{1}+B^{1}C_{g}^{2}-B^{2}C_{f}^{2}+B^{2}C_{g}^{1}-
\notag \\
2B^{1}B_{x}^{1}-2B^{1}B_{y}^{2}+2B^{2}B_{x}^{2}-2B^{2}B_{y}^{1}-2B_{xf}^{1}-2B_{yf}^{2}-2B_{xg}^{2}+2B_{yg}^{1}=0,
\tag{29}
\end{gather}%
\begin{gather}
3A_{xx}^{2}-3A_{yy}^{2}-6A_{xy}^{1}+3C^{2}A_{x}^{1}+3C^{2}A_{y}^{2}+3A_{x}^{2}C^{1}-3C^{1}A_{y}^{1}-3D^{2}A_{f}^{1}-3D^{2}A_{g}^{2}-
\notag \\
3D^{1}A_{f}^{2}+3D^{1}A_{g}^{1}+3A^{2}C_{x}^{1}+3A^{2}C_{y}^{2}+3A^{1}C_{y}^{2}-3A^{1}C_{x}^{1}+C_{ff}^{2}-C_{gg}^{2}-2C_{fg}^{1}-
\notag \\
6A^{2}D_{f}^{1}-6A^{2}D_{g}^{2}-6A^{1}D_{f}^{2}+6A^{1}D_{g}^{1}+B^{2}C_{f}^{1}+B^{2}C_{g}^{2}+B^{1}C_{f}^{2}-B^{1}C_{g}^{1}-
\notag \\
2B^{2}B_{x}^{1}-2B^{2}B_{y}^{2}-2B^{1}B_{x}^{2}+2B^{1}B_{y}^{1}-2B_{xf}^{2}+2B_{yf}^{1}+2B_{xg}^{1}-2B_{yg}^{2}=0
\tag{30}
\end{gather}%
\begin{gather}
6D^{1}A_{x}^{1}+6D^{1}A_{y}^{2}-6D^{2}A_{x}^{2}+6D^{2}A_{y}^{1}-3D^{1}B_{f}^{1}-3D^{1}B_{g}^{2}+3D^{2}B_{f}^{2}-3D^{2}B_{g}^{1}+
\notag \\
3A^{1}D_{x}^{1}+3A^{1}D_{y}^{2}-3A^{2}D_{x}^{2}+3A^{2}D_{y}^{1}+B_{xx}^{1}-B_{yy}^{1}+2B_{xy}^{2}-2C_{xf}^{1}-2C_{yf}^{2}-
\notag \\
2C_{xg}^{2}+2C_{yg}^{1}-3B^{1}D_{f}^{1}-3B^{1}D_{g}^{2}+3B^{2}D_{g}^{2}-3B^{2}D_{g}^{1}+3D_{ff}^{1}-3D_{gg}^{1}+6D_{fg}^{2}+
\notag \\
2C^{1}C_{f}^{1}+2C^{1}C_{g}^{2}-2C^{2}C_{f}^{2}+2C^{2}C_{g}^{1}-C^{1}B_{x}^{1}-C^{1}B_{y}^{2}+C^{2}B_{x}^{2}-C^{2}B_{y}^{1}=0,
\tag{31}
\end{gather}%
\begin{gather}
6D^{2}A_{x}^{1}+6D^{2}A_{y}^{2}+6D^{1}A_{x}^{2}-6D^{1}A_{y}^{1}-3D^{2}B_{f}^{1}-3D^{2}B_{g}^{2}-3D^{1}B_{f}^{2}+3D^{1}B_{g}^{1}+
\notag \\
3A^{2}D_{x}^{1}+3A^{2}D_{y}^{2}+3A^{1}D_{x}^{2}-3A^{1}D_{y}^{1}+B_{xx}^{2}-B_{yy}^{2}-2B_{xy}^{1}-2C_{xf}^{2}+2C_{yf}^{1}+
\notag \\
2C_{xg}^{1}+2C_{yg}^{2}-3B^{2}D_{f}^{1}-3B^{2}D_{g}^{2}-3B^{1}D_{f}^{2}+3B^{1}D_{g}^{1}+3D_{ff}^{2}-3D_{gg}^{2}-6D_{fg}^{1}+
\notag \\
2C^{2}C_{f}^{1}-2C^{2}C_{g}^{2}+2C^{1}C_{f}^{2}-2C^{1}C_{g}^{1}-C^{2}B_{x}^{1}-C^{2}B_{y}^{2}-C^{1}B_{x}^{2}+C^{1}B_{y}^{1}=0;
\tag{32}
\end{gather}%
\textbf{3.} The system of PDEs corresponding to a CODE has four real
symmetries $\mathbf{X}_{1},~\mathbf{Y}_{1},~\mathbf{X}_{2}$ and $\mathbf{Y}%
_{2}$ with 
\begin{equation}
\mathbf{X}_{1}=\rho _{1}\mathbf{X}_{2}-\rho _{2}\mathbf{Y}_{2},~\mathbf{Y}%
_{1}\ =\rho _{1}\mathbf{Y}_{2}+\rho _{2}\mathbf{X}_{2},  \tag{33}
\end{equation}%
for nonconstant $\rho _{1}$ and $\rho _{2}$ and they satisfy%
\begin{equation}
\lbrack \mathbf{X}_{1},\mathbf{X}_{2}]-[\mathbf{Y}_{1},\mathbf{Y}_{2}]=0,~[%
\mathbf{X}_{1},\mathbf{Y}_{2}]+[\mathbf{Y}_{1},\mathbf{X}_{2}]=0,  \tag{34}
\end{equation}%
such that a point transformation $(x,y,f,g)\longrightarrow (X,Y,F,G),$ which
brings $\mathbf{X}_{1},~\mathbf{Y}_{1},~\mathbf{X}_{2}$ and $\mathbf{Y}_{2}$
to their canonical form%
\begin{equation}
\mathbf{X}_{1}=\frac{\partial }{\partial F},~\mathbf{Y}_{1}=\frac{\partial }{%
\partial G},~\mathbf{X}_{2}=X\frac{\partial }{\partial F}+Y\frac{\partial }{%
\partial G},~\mathbf{Y}_{2}=Y\frac{\partial }{\partial F}-X\frac{\partial }{%
\partial G}  \tag{35}
\end{equation}%
reduces the system (27) and (28) to the linear form%
\begin{equation}
F_{XX}-F_{YY}+2G_{XY}=W_{1}(X,Y),  \tag{36}
\end{equation}%
\begin{equation}
G_{XX}-G_{YY}-2F_{XY}=W_{2}(X,Y);  \tag{37}
\end{equation}%
\textbf{4.} The system of PDEs corresponding to a CODE has four real
symmetries $\mathbf{X}_{1},~\mathbf{Y}_{1},~\mathbf{X}_{2}$ and $\mathbf{Y}%
_{2}$ with 
\begin{equation}
\mathbf{X}_{1}=\rho _{1}\mathbf{X}_{2}-\rho _{2}\mathbf{Y}_{2},~\mathbf{Y}%
_{1}\ =\rho _{1}\mathbf{Y}_{2}+\rho _{2}\mathbf{X}_{2},  \tag{38}
\end{equation}%
for nonconstant $\rho _{1}$ and $\rho _{2}$ and they satisfy either 
\begin{equation}
\lbrack \mathbf{X}_{1},\mathbf{X}_{2}]-[\mathbf{Y}_{1},\mathbf{Y}_{2}]\neq 0%
\text{ or }[\mathbf{X}_{1},\mathbf{Y}_{2}]+[\mathbf{Y}_{1},\mathbf{X}%
_{2}]\neq 0  \tag{39}
\end{equation}%
such that a point transformation $(x,y,f,g)\longrightarrow (X,Y,F,G),$ which
brings $\mathbf{X}_{1},~\mathbf{Y}_{1},~\mathbf{X}_{2}$ and $\mathbf{Y}_{2}$
to their canonical form%
\begin{equation}
\mathbf{X}_{1}=\frac{\partial }{\partial F},~\mathbf{Y}_{1}=\frac{\partial }{%
\partial G},~\mathbf{X}_{2}=F\frac{\partial }{\partial F}+G\frac{\partial }{%
\partial G},~\mathbf{Y}_{2}=G\frac{\partial }{\partial F}-F\frac{\partial }{%
\partial G},  \tag{40}
\end{equation}%
reduces the system (27) and (28) to the linear form%
\begin{equation}
F_{XX}-F_{YY}+2G_{XY}=HW_{1}(X,Y)-LW_{2}(X,Y),  \tag{41}
\end{equation}%
\begin{equation}
G_{XX}-G_{YY}-2F_{XY}=HW_{2}(X,Y).+LW_{1}(X,Y).  \tag{42}
\end{equation}%
where 
\begin{equation}
H=F_{X}+G_{Y},~~L=G_{X}-F_{Y};  \tag{43}
\end{equation}%
\textbf{5.} The system of PDEs corresponding to a CODE has four real
symmetries $\mathbf{X}_{1},~\mathbf{Y}_{1},~\mathbf{X}_{2}$ and $\mathbf{Y}%
_{2}$ with 
\begin{equation}
\mathbf{X}_{1}\neq \rho _{1}\mathbf{X}_{2}-\rho _{2}\mathbf{Y}_{2},~\mathbf{Y%
}_{1}\ \neq \rho _{1}\mathbf{Y}_{2}+\rho _{2}\mathbf{X}_{2},  \tag{44}
\end{equation}%
for nonconstant $\rho _{1}$ and $\rho _{2}$ and they satisfy%
\begin{equation}
\lbrack \mathbf{X}_{1},\mathbf{X}_{2}]-[\mathbf{Y}_{1},\mathbf{Y}_{2}]=0,~[%
\mathbf{X}_{1},\mathbf{Y}_{2}]+[\mathbf{Y}_{1},\mathbf{X}_{2}]=0,  \tag{45}
\end{equation}%
such that a point transformation $(x,y,f,g)\longrightarrow (X,Y,F,G),$ which
brings $\mathbf{X}_{1},~\mathbf{Y}_{1},~\mathbf{X}_{2}$ and $\mathbf{Y}_{2}$
to their canonical form%
\begin{equation}
\mathbf{X}_{1}=\frac{\partial }{\partial X},~\mathbf{Y}_{1}=\frac{\partial }{%
\partial Y},~\mathbf{X}_{2}=\frac{\partial }{\partial F},~\mathbf{Y}_{2}=%
\frac{\partial }{\partial G}  \tag{46}
\end{equation}%
reduces the system (27) and (28) to the system of PDEs corresponding to CODE
is at most cubic in all its first derivatives;\newline
\textbf{6.} The system of PDEs corresponding to a CODE has four real
symmetries $\mathbf{X}_{1},~\mathbf{Y}_{1},~\mathbf{X}_{2}$ and $\mathbf{Y}%
_{2}$ with 
\begin{equation}
\mathbf{X}_{1}\neq \rho _{1}\mathbf{X}_{2}-\rho _{2}\mathbf{Y}_{2},~\mathbf{Y%
}_{1}\ \neq \rho _{1}\mathbf{Y}_{2}+\rho _{2}\mathbf{X}_{2},  \tag{47}
\end{equation}%
for nonconstant $\rho _{1}$ and $\rho _{2}$ and they satisfy either 
\begin{equation}
\lbrack \mathbf{X}_{1},\mathbf{X}_{2}]-[\mathbf{Y}_{1},\mathbf{Y}_{2}]\neq 0%
\text{ or }[\mathbf{X}_{1},\mathbf{Y}_{2}]+[\mathbf{Y}_{1},\mathbf{X}%
_{2}]\neq 0  \tag{48}
\end{equation}%
such that a point transformation $(x,y,f,g)\longrightarrow (X,Y,F,G),$ which
brings $\mathbf{X}_{1},~\mathbf{Y}_{1},~\mathbf{X}_{2}$ and $\mathbf{Y}_{2}$
to their canonical form%
\begin{equation}
\mathbf{X}_{1}=\frac{\partial }{\partial F},~\mathbf{Y}_{1}=\frac{\partial }{%
\partial G},~\mathbf{X}_{2}=X\frac{\partial }{\partial X}+Y\frac{\partial }{%
\partial Y}+F\frac{\partial }{\partial F}+G\frac{\partial }{\partial G},{} 
\tag{49}
\end{equation}%
\begin{equation}
~\mathbf{Y}_{2}=Y\frac{\partial }{\partial X}-X\frac{\partial }{\partial Y}+G%
\frac{\partial }{\partial F}-F\frac{\partial }{\partial G},  \tag{50}
\end{equation}%
reduces the system (27) and (28) to the linear form%
\begin{gather}
X(F_{XX}-F_{YY}+2G_{XY})-Y(G_{XX}-G_{YY}-2F_{XY})=a_{1}(H^{3}-3HL^{2})- 
\notag \\
a_{2}(3H^{2}L-L^{3})+b_{1}(H^{2}-L^{2})-2b_{2}HL+\frac{1}{%
3(a_{1}^{2}+a_{2}^{2})}\{3(a_{1}^{2}+a_{2}^{2})+(b_{1}^{2}-b_{2}^{2})a_{1}+ 
\notag \\
2b_{1}b_{2}a_{2}\}H-\frac{1}{3(a_{1}^{2}+a_{2}^{2})}%
\{2b_{1}b_{2}a_{1}-a_{2}(b_{1}^{2}-b_{2}^{2})\}L+\frac{1}{%
3(a_{1}^{2}+a_{2}^{2})}(b_{1}a_{1}+b_{2}a_{2})+  \notag \\
\frac{1}{27(a_{1}^{2}+a_{2}^{2})}%
\{(b_{1}^{3}-3b_{1}b_{2}^{2})(a_{1}^{2}-a_{2}^{2})+2a_{1}a_{2}(3b_{1}^{2}b_{2}-b_{2}^{3}),
\tag{51}
\end{gather}%
\begin{gather}
Y(F_{XX}-F_{YY}+2G_{XY})+X(G_{XX}-G_{YY}-2F_{XY})=a_{2}(H^{3}-3HL^{2})+ 
\notag \\
a_{1}(3H^{2}L-L^{3})+b_{2}(H^{2}-L^{2})+2b_{1}HL+\frac{1}{%
3(a_{1}^{2}+a_{2}^{2})}\{3(a_{1}^{2}+a_{2}^{2})+(b_{1}^{2}-b_{2}^{2})a_{1}+ 
\notag \\
2b_{1}b_{2}a_{2}\}L+\frac{1}{3(a_{1}^{2}+a_{2}^{2})}%
\{2b_{1}b_{2}a_{1}-a_{2}(b_{1}^{2}-b_{2}^{2})\}H+\frac{1}{%
3(a_{1}^{2}+a_{2}^{2})}(b_{2}a_{1}-b_{1}a_{2})+  \notag \\
\frac{1}{27(a_{1}^{2}+a_{2}^{2})}%
\{(3b_{1}^{2}b_{2}-b_{2}^{3})(a_{1}^{2}-a_{2}^{2})-2(b_{1}^{3}-3b_{1}b_{2}^{2})a_{1}a_{2}\}.
\tag{52}
\end{gather}%
The invertible real transformations%
\begin{equation}
\tilde{X}=F+\frac{1}{3(a_{1}^{2}+a_{2}^{2})}%
\{(b_{1}a_{1}+b_{2}a_{2})X-(b_{2}a_{1}-b_{1}a_{2})Y\},  \tag{53}
\end{equation}%
\begin{equation}
~\tilde{Y}=G+\frac{1}{3(a_{1}^{2}+a_{2}^{2})}%
\{(b_{1}a_{1}+b_{2}a_{2})Y+(b_{2}a_{1}-b_{1}a_{2})X\},  \tag{54}
\end{equation}%
\begin{equation}
\tilde{F}=\frac{1}{2}(F^{2}-G^{2})+\frac{1}{3(a_{1}^{2}+a_{2}^{2})}%
[\{(b_{1}a_{1}+b_{2}a_{2})X-(b_{2}a_{1}-b_{1}a_{2})Y\}F-  \notag
\end{equation}%
\begin{equation}
\{(b_{1}a_{1}+b_{2}a_{2})Y+(b_{2}a_{1}-b_{1}a_{2})X\}G]+\frac{1}{%
18(a_{1}^{2}+a_{2}^{2})^{2}}[\{(b_{1}^{2}-b_{2}^{2})(a_{1}^{2}-a_{2}^{2})+ 
\notag
\end{equation}%
\begin{equation}
4b_{1}b_{2}a_{1}a_{2}\}(X^{2}-Y^{2})-2XY%
\{(2b_{1}b_{2}(a_{1}^{2}-a_{2}^{2})-2a_{1}a_{2}(b_{1}^{2}-b_{2}^{2})\}]+ 
\notag
\end{equation}%
\begin{equation}
\frac{1}{2(a_{1}^{2}+a_{2}^{2})}\{a_{1}(X^{2}-Y^{2})+2a_{2}XY\},  \tag{55}
\end{equation}%
\begin{equation}
\tilde{G}=FG+\frac{1}{3(a_{1}^{2}+a_{2}^{2})}%
[\{(b_{1}a_{1}+b_{2}a_{2})X-(b_{2}a_{1}-b_{1}a_{2})Y\}G+  \notag
\end{equation}%
\begin{equation}
\{(b_{1}a_{1}+b_{2}a_{2})Y+(b_{2}a_{1}-b_{1}a_{2})X\}F]+\frac{1}{%
18(a_{1}^{2}+a_{2}^{2})^{2}}[2XY\{(b_{1}^{2}-b_{2}^{2})(a_{1}^{2}-a_{2}^{2})+
\notag
\end{equation}%
\begin{equation}
4b_{1}b_{2}a_{1}a_{2}\}-(X^{2}-Y^{2})%
\{2b_{1}b_{2}(a_{1}^{2}-a_{2}^{2})-2a_{1}a_{2}(b_{1}^{2}-b_{2}^{2})\}]+ 
\notag
\end{equation}%
\begin{equation}
\frac{1}{2(a_{1}^{2}+a_{2}^{2})}\{2XYa_{1}-a_{2}(X^{2}-Y^{2})\}.  \tag{56}
\end{equation}%
transform (51) and(52) into the system of linear PDEs 
\begin{equation}
\tilde{F}_{\tilde{X}~\tilde{X}}-\tilde{F}_{\tilde{Y}~\tilde{Y}}+2~\tilde{G}_{%
\tilde{X}~\tilde{Y}}=0,  \tag{57}
\end{equation}%
\begin{equation}
\tilde{G}_{\tilde{X}~\tilde{X}}-\tilde{G}_{\tilde{Y}~\tilde{Y}}-2~\tilde{F}_{%
\tilde{X}~\tilde{Y}}=0.  \tag{58}
\end{equation}%
\newline
\newline
\textbf{Examples.} Now we discuss some illustrative examples.\newline
\textbf{1.} Consider a nonlinear CODE of the form%
\begin{equation}
u^{\prime \prime }+3u~u^{\prime }+u^{3}=0,  \tag{59}
\end{equation}%
which is linearizable as it satisfies the Lie complex conditions. The set of
two noncommutating CLSs is%
\begin{equation}
\mathbf{Z}_{1}=\frac{\partial }{\partial z},~\mathbf{Z}_{2}=z\frac{\partial 
}{\partial z}-u\frac{\partial }{\partial u}.  \tag{60}
\end{equation}%
Note that $\mathbf{Z}_{1}\neq \rho (z,u)\mathbf{Z}_{2}.$ We invoke condition
nine of Theorem 1 to find a linearization transformation. The complex point
transformation that reduces the symmetries (60) to their canonical form is 
\begin{equation}
Z=\frac{1}{u},~U=z+\frac{1}{u}  \tag{61}
\end{equation}%
and (59) reduces to 
\begin{equation}
ZU^{\prime \prime }=-U^{\prime 3}+6U^{\prime 2}-11U^{\prime }+6  \tag{62}
\end{equation}%
by means of the transformation (61). Equation (59) linearizes to $\tilde{U}%
^{\prime \prime }=0$ via the complex transformations (25) with $a=-1,~b=6$
and by using (60). That is 
\begin{equation}
\tilde{Z}=z-\frac{1}{u},~\tilde{U}=\frac{z^{2}}{2}-\frac{z}{u}.  \tag{63}
\end{equation}%
It may be seen that the above transformation is a proper complex
transformation thus one does not get a contradiction as was in P - I \cite%
{Ali2}. Further, the procedure that seems to be analogous to analytic
continuation do not arise here. The system of PDEs corresponding to (59) is 
\begin{equation}
f_{xx}-f_{yy}+2g_{xy}=-3(fh-gl)-(f^{3}-3fg^{2}),  \tag{64}
\end{equation}%
\begin{equation}
g_{xx}-g_{yy}-2f_{xy}=-3(gh+fl)-(3f^{2}g-g^{3}),  \tag{65}
\end{equation}%
which admits the RLSs 
\begin{equation}
\mathbf{X}_{1} =\frac{\partial }{\partial x},\mathbf{Y}_{1}=-\frac{%
\partial }{\partial y},  \tag{66} \\
\end{equation}%
\begin{equation}
\mathbf{X}_{2} =x\frac{\partial }{\partial x}+y\frac{\partial }{\partial y}%
-f\frac{\partial }{\partial f}-g\frac{\partial }{\partial g},  \tag{67} \\
\end{equation}%
\begin{equation}
\mathbf{Y}_{2} =y\frac{\partial }{\partial x}-x\frac{\partial }{\partial y}%
-g\frac{\partial }{\partial f}+f\frac{\partial }{\partial g}.  \tag{68}
\end{equation}%
These generators satisfy (47) and (49). The transformation%
\begin{equation}
X =\frac{f}{f^{2}+g^{2}},\text{ \ \ }Y=\frac{-g}{f^{2}+g^{2}},  \tag{69}
\\
\end{equation}%
\begin{equation}
F =x+\frac{f}{f^{2}+g^{2}},G=y-\frac{g}{f^{2}+g^{2}},  \tag{70}
\end{equation}%
transforms (64) and (65) into%
\begin{gather}
X(F_{XX}-F_{YY}+2G_{XY})-Y(G_{XX}-G_{YY}-2F_{XY})=-H^{3}+3HL^{2}-  \notag \\
\text{ \ \ \ \ \ \ \ \ \ \ \ \ \ \ \ \ \ \ \ \ \ \ \ \ \ \ \ \ \ \ \ \ \ \ \
\ \ \ \ \ \ \ \ \ \ \ \ \ \ \ \ \ \ \ \ \ \ \ \ \ \ \ \ \ \ \ \ }6H^{2}-11H+6
\tag{71} \\
X(G_{XX}-G_{YY}-2F_{XY})+Y(F_{XX}-F_{YY}+2G_{XY})=L^{3}-3H^{2}L+  \notag \\
\text{ \ \ \ \ \ \ \ \ \ \ \ \ \ \ \ \ \ \ \ \ \ \ \ \ \ \ \ \ \ \ \ \ \ \ \
\ \ \ \ \ \ \ \ \ \ \ \ \ \ \ \ \ \ \ \ \ \ \ \ \ \ \ }12HL-11L  \tag{72}
\end{gather}%
where $H=F_{X}+G_{Y},~L=G_{X}-F_{Y}.$ The above system can further be
reduced into very simplified linear PDEs by using 
\begin{equation}
\tilde{X} =x-\frac{f}{f^{2}+g^{2}},\text{ }\tilde{Y}=y+\frac{g}{f^{2}+g^{2}%
},  \tag{73} \\
\end{equation}%
\begin{equation}
\tilde{F} =\frac{1}{2}(x^{2}-y^{2})-\frac{1}{f^{2}+g^{2}}(xf+yg), 
\tag{74} \\
\end{equation}%
\begin{equation}
\tilde{G} =xy-\frac{1}{f^{2}+g^{2}}(yf-xg).  \tag{75}
\end{equation}%
The linear PDEs are 
\begin{equation}
\tilde{F}_{\tilde{X}~\tilde{X}}-\tilde{F}_{\tilde{Y}~\tilde{Y}}+2~\tilde{G}_{%
\tilde{X}~\tilde{Y}} =0,  \tag{76} \\
\end{equation}%
\begin{equation}
\tilde{G}_{\tilde{X}~\tilde{X}}-\tilde{G}_{\tilde{Y}~\tilde{Y}}-2~\tilde{F}_{%
\tilde{X}~\tilde{Y}} =0.  \tag{77}
\end{equation}%
\newline
\newline
\textbf{2.} Consider a second order nonlinear CODE with an arbitrary
function $w(z)$%
\begin{equation}
u~u^{\prime \prime }=u^{\prime 2}+w(z)u^{2}.  \tag{78}
\end{equation}%
The above CODE admits two CLSs of the form%
\begin{equation}
\mathbf{Z}_{1}=zu\frac{\partial }{\partial u},~\mathbf{Z}_{2}=u\frac{%
\partial }{\partial u},  \tag{79}
\end{equation}%
and thus%
\begin{equation}
\lbrack \mathbf{Z}_{1},\mathbf{Z}_{2}]=0\text{ and }\mathbf{Z}_{2}=\frac{1}{z%
}\mathbf{Z}_{1}.  \tag{80}
\end{equation}%
By using the complex transformation%
\begin{equation}
Z=\frac{1}{z}\text{ and }U=\frac{1}{z}\log u  \tag{81}
\end{equation}%
(78) can be reduced into a linear CODE%
\begin{equation}
U^{\prime \prime }=\frac{1}{Z^{3}}w(\frac{1}{Z}).  \tag{82}
\end{equation}%
The system of PDEs corresponding to (78) are 
\begin{equation}
f(f_{xx}{\small -}f_{yy}{\small +}2g_{xy}){\small -}g(g_{xx}{\small -}g_{yy}%
{\small -}2f_{xy}) {\small =}h^{2}{\small -}l^{2}{\small +}w_{1}(f^{2}%
{\small -}g^{2}){\small -}2fgw_{2},  \tag{83} \\
\end{equation}%
\begin{equation}
f(g_{xx}{\small -}g_{yy}{\small -}2f_{xy}){\small +}(f_{xx}{\small -}f_{yy}%
{\small +}2g_{xy})g {\small =}2hl{\small +}2w_{1}fg{\small +}w_{2}(f^{2}%
{\small -}g^{2}),  \tag{84}
\end{equation}%
where%
\begin{equation*}
h=f_{x}+g_{y},l=g_{x}-f_{y}.
\end{equation*}%
Invoking the transformation 
\begin{equation}
X =\frac{x}{x^{2}+y^{2}}\text{, }Y=\frac{-y}{x^{2}+y^{2}}  \tag{85} \\
\end{equation}%
\begin{equation}
F =\frac{1}{x^{2}+y^{2}}[\frac{1}{2}x\ln (f^{2}+g^{2})+y\tan ^{-1}(\frac{g%
}{f})],  \tag{86} \\
\end{equation}%
\begin{equation}
G =\frac{1}{x^{2}+y^{2}}[x\tan ^{-1}(\frac{g}{f})-\frac{y}{2}\ln
(f^{2}+g^{2}).  \tag{87}
\end{equation}%
Equations (83) and (84) reduce to the linear system of PDEs 
\begin{equation}
F_{XX}{\small -}F_{YY}{\small +}2G_{XY} {\small =}\frac{1}{X^{2}+Y^{2}}%
[(X^{3}{\small -}3XY^{2})w_{1}{\small -}(Y^{3}{\small -}3X^{2}Y)w_{2}], 
\tag{88} \\
\end{equation}%
\begin{equation}
G_{XX}{\small -}G_{YY}{\small -}2F_{XY} {\small =}\frac{1}{X^{2}+Y^{2}}%
[(X^{3}{\small -}3XY^{2})w_{2}{\small +}(Y^{3}{\small -}3X^{2}Y)w_{1}], 
\tag{89}
\end{equation}%
where%
\begin{equation}
w_{1}=w_{1}(\frac{X}{X^{2}+Y^{2}},\frac{-Y}{X^{2}+Y^{2}}),\text{ }%
w_{2}=w_{2}(\frac{X}{X^{2}+Y^{2}},\frac{-Y}{X^{2}+Y^{2}}).  \tag{90}
\end{equation}%
\newline
\newline
\textbf{3. }Consider the nonlinear CODE%
\begin{equation}
u^{\prime \prime }=\frac{1}{z}(u^{\prime }+u^{^{\prime }3}).  \tag{91}
\end{equation}%
This equation has CLSs%
\begin{equation}
\mathbf{Z}_{1}=\frac{1}{z}\frac{\partial }{\partial z},~\mathbf{Z}_{2}=\frac{%
u}{z}\frac{\partial }{\partial z}.  \tag{92}
\end{equation}%
The complex transformation%
\begin{equation}
U=\frac{1}{2}z^{2}\text{ and }Z=u,  \tag{93}
\end{equation}%
transforms (91) into a linear CODE%
\begin{equation}
U^{\prime \prime }+1=0.  \tag{94}
\end{equation}%
The system of PDEs corresponding to (91) is 
\begin{equation}
f_{xx}-f_{yy}+2g_{xy} =\frac{1}{x^{2}+y^{2}}%
[xh(1+h^{2}-3l^{2})+yl(1-l^{2}+3h^{2})],  \tag{95} \\
\end{equation}%
\begin{equation}
g_{xx}-g_{yy}-2f_{xy} =\frac{1}{x^{2}+y^{2}}%
[xl(1-l^{2}+3h^{2})+yh(1+h^{2}-3l^{2})],  \tag{96}
\end{equation}%
and admits CLSs%
\begin{equation}
\mathbf{X}_{1} =\frac{x}{x^{2}+y^{2}}\frac{\partial }{\partial x}-\frac{y}{%
x^{2}+y^{2}}\frac{\partial }{\partial y},\mathbf{Y}_{1}=\frac{-y}{x^{2}+y^{2}%
}\frac{\partial }{\partial x}-\frac{x}{x^{2}+y^{2}}\frac{\partial }{\partial
y},  \tag{97} \\
\end{equation}%
\begin{equation}
\mathbf{X}_{2} =\frac{fx}{x^{2}+y^{2}}\frac{\partial }{\partial x}-\frac{fy%
}{x^{2}+y^{2}}\frac{\partial }{\partial y}+\frac{yg}{x^{2}+y^{2}}\frac{%
\partial }{\partial x}-\frac{gx}{x^{2}+y^{2}}\frac{\partial }{\partial y}, 
\tag{98} \\
\end{equation}%
\begin{equation}
\mathbf{Y}_{2} =\frac{gx}{x^{2}+y^{2}}\frac{\partial }{\partial x}-\frac{gy%
}{x^{2}+y^{2}}\frac{\partial }{\partial y}-\frac{fy}{x^{2}+y^{2}}\frac{%
\partial }{\partial x}-\frac{fx}{x^{2}+y^{2}}\frac{\partial }{\partial y}. 
\tag{99}
\end{equation}%
The above RLSs correspond to the condtion five of Theorm 2. The real
transformation that linearize system (95) and (96) is%
\begin{equation}
X=f,Y=g,\text{ \ }F=\frac{1}{2}(x^{2}-y^{2}),G=xy.  \tag{100}
\end{equation}%
The linearized system of PDEs is 
\begin{gather}
F_{XX}-F_{YY}+2G_{XY}+1=0,  \tag{101} \\
G_{XX}-G_{YY}-2F_{XY}=0.  \tag{102}
\end{gather}%
\newline
\newline
\textbf{4. }Consider the nonlinear CODE%
\begin{equation}
u^{\prime \prime }=1+(u^{\prime }-z)^{2}w(2u-z^{2})  \tag{103}
\end{equation}%
which admits the CLSs%
\begin{equation}
\mathbf{Z}_{1}=\frac{\partial }{\partial z}+z\frac{\partial }{\partial u},~%
\mathbf{Z}_{2}=z\frac{\partial }{\partial z}+z^{2}\frac{\partial }{\partial u%
}.  \tag{104}
\end{equation}%
Equation (103) reduces to a linear CODE by using the complex transformation%
\begin{equation}
Z=2u-z^{2}\text{ and }U=z.  \tag{105}
\end{equation}%
It becomes%
\begin{equation}
U^{\prime \prime }=\frac{-1}{2}U^{\prime }w(Z).  \tag{106}
\end{equation}%
The system of PDEs corresponding to (103) is 
\begin{gather}
f_{xx}-f_{yy}+2g_{xy}\text{ }{\small =}\text{ }1+\{(h-x)^{2}-(l-y)^{2}%
\}w_{1}-2(h-x)(l-y)w_{2},  \tag{107} \\
g_{xx}-g_{yy}-2f_{xy}\text{ }{\small =}\text{ }2(h-x)(l-y)w_{1}+%
\{(h-x)^{2}-(l-y)^{2}\}w_{2},  \tag{108}
\end{gather}%
which admits RLSs%
\begin{equation}
\mathbf{X}_{1} =\frac{\partial }{\partial x}+x\frac{\partial }{\partial f}%
+y\frac{\partial }{\partial g},\mathbf{Y}_{1}=x\frac{\partial }{\partial g}-y%
\frac{\partial }{\partial f}-\frac{\partial }{\partial y},  \tag{109} \\
\end{equation}%
\begin{equation}
\mathbf{X}_{2} =x\frac{\partial }{\partial x}+y\frac{\partial }{\partial y}%
+(x^{2}-y^{2})\frac{\partial }{\partial f}+2xy\frac{\partial }{\partial g}, 
\tag{110} \\
\end{equation}%
\begin{equation}
\mathbf{Y}_{2} =y\frac{\partial }{\partial x}-x\frac{\partial }{\partial y}%
+2xy\frac{\partial }{\partial f}-(x^{2}-y^{2})\frac{\partial }{\partial g}, 
\tag{111}
\end{equation}%
which satisfy condition six of Theorm 2. The transformation%
\begin{equation}
X =2f-x^{2}+y^{2},~Y=2g-2xy,  \tag{112} \\
\end{equation}%
\begin{equation}
F =x,\text{ \ }G=y.  \tag{113}
\end{equation}%
yields system of linear PDEs%
\begin{equation}
F_{XX}-F_{YY}+2G_{XY} =-\frac{1}{2}(Hw_{1}-Lw_{2}),  \tag{114} \\
\end{equation}%
\begin{equation}
G_{XX}-G_{YY}-2F_{XY} =-\frac{1}{2}(Hw_{2}+Lw_{1}),  \tag{115}
\end{equation}%
where 
\begin{equation}
w(z)=w_{1}+iw_{2},~\ H=F_{X}+G_{Y},\text{ \ }L=G_{X}-F_{Y}.  \tag{116}
\end{equation}

\section{\textbf{Conclusion}}

\qquad Finding exact solutions of nonlinear DEs is difficult. Linearization
of DEs is a method to convert the original nonlinear DE by an equivalent
linear DE so that exact solutions can be constructed directly. But it
requires the existence and construction of such transformations (point,
tangent, contact or complex) that transforms nonlinear DEs into linear DEs.
For scalar ODEs Lie provided nontrivial ways of constructing point
transformations via symmetries. But it is difficult to extend the Lie
conditions for systems of ODEs and for PDEs (see e.g., \cite{blu, blu1, kum}%
, $[15-19$]). Finding transformations that map systems of nonlinear ODEs and
PDEs into linear systems is not only tedious but also highly nontrivial.

We have applied CLS analysis to find non-trivial ways of \textit{reducing
order, linearizing and solving }certain systems of ODEs and PDEs. Several \
known results of classical symmetry method were studied in the complex
domain. Various results for systems of PDEs were obtained, e.g. the CLS
analysis gets used to study CLSs for CODEs and respective real Lie
symmetries for systems of PDEs corresponding to CODEs is discussed in \cite%
{Ali}. Variational problems i.e. complex Lagrangians of CODEs results in
Lagrangians of systems of PDEs corresponding to those CODEs and double
reduction in order of CODEs via complex Noether symmetries are studied in 
\cite{Ali1}. Further, variational problems for r-CODEs and their use in
constructing Lagrangians for systems of ODEs are discussed in \cite{Ali3}.

In this work we have looked at linearizability criteria for second order
CODEs. We obtained analogous Lie conditions for CODEs and applied them to
the linearization of systems of nonlinear PDEs. The corresponding statement
for a system of two PDEs (in fact four i.e. by incorporating linear
CR-equations) associated to a second order CODE was presented. Examples were
given for the linearization of systems of PDEs including the construction of
the point linearzing transformations.

It is interesting to see that in this paper the complex function is a
function of complex plane that requires CR-equations to hold for complex
differentiability. Thus, we mapped complex solutions of CODEs into complex
solutions via complex transformations. Since, CR-equations are linear thus
makes no difference in the linearization of systems of PDEs corresponding to
CODEs except that all the symmetries are analytic i.e. their coefficents
satisfy CR-equations. These form the conformal subalgebra admitted by
systems of PDEs. It is hoped that further classifications of systems of PDEs
in \cite{Ali} with respect to conformal algebras admitted by them can be
done. 

We have seen that linearization theorems can be carried over from the real
case to the complex case with remarkable results for the linearization of
systems of PDEs. Thus we were able to linearize those systems of nonlinear
PDEs that correspond to some CODEs. There are several possible ways of using
CLS analysis. It would be of great interest to extend the geometric proof of
linearization for RODEs \cite{ibr2, mah2} to CODEs. The analysis can also be
used for classification of those systems of nonlinear PDEs that correpsond
to CODEs. An extension of the Lie linearization results for systems of
quadratically and cubically semilinear CODEs in \cite{mah, mah2}, in the
complex domain may yield useful and highly nontrivial results for systems of
PDEs. Further, results of linearization for third and fourth order CODEs may
be obtained by using results for third and fourth order RODEs in \cite{mah4,
mah3, mah5}$.$

It is suggested that the Lie's criteria for linearization of two CODEs can
be found in a similar way, which can then be used to linearize four
nonlinear PDEs corresponding to a system of two CODEs. A system of two
second order RODES admits $5,6,7,8$ or $15$ RLSs and the maximal symmetry
algebra is $sl(4,R)$ for the simplest system \cite{waf}. Thus, it may be
argued that the system of two second order CODEs also admit $5,6,7,8$ or $15$
CLSs. These result in $10,12,14,16$ or $30$ RLSs which form the subalgebra
of symmetries admitted by four second order PDEs that correspond to two
second order CODEs. Further, maximal symmetry algebra $sl(4,C)$ is admitted
by the simplest systems of two second order CODEs.

It is also hoped that \textquotedblleft complex magic\textquotedblright\ can
be further extended to \textquotedblleft \textit{hypercomplex magic}%
\textquotedblright\ by introducing hypercomplex variables (e.g. quaternions,
octonions, Clifford or Grassmanian variables) that may result in
hypercomplex Lie symmetries. Then, these can be used to reduce the order,
linearize and solve large classes of systems of PDEs. It is hoped that
non-commutative behavior of hypercomplex variables may yield non-trivial
generalization of Lie's work. Further, use of the hypercomplex Lie
symmetries in geometric calculus \cite{hes} may result in non-trivial
physical implications.

In paper - I, we have also used CLS analysis in constructing linearizability
criteria for certain systems of nonlinear RODEs by considering r-CODEs. An\
r-CODE\ is complex ordinary differential equations in a single complex
function of one real variable that yields a system of two nonlinear RODEs.
The linearization of r-CODEs directly gives linearizability criteria for
systems of RODEs corresonding to those r-CODEs \cite{Ali2}.

\section{\textbf{Acknowledgments}}

\qquad SA is most grateful to NUST and DECMA in providing finnacial
assistance for his stay at Wits University, South Africa, where this work
was done. AQ acknowledges the School of Computational and Applied
Mathematics (DECMA) for funding his stay at the university.

\end{document}